\input amstex 
\documentstyle{amsppt} 
\loadbold
\magnification=1200

\pagewidth{6.43truein}
\pageheight{8.5truein}

\define\C{{\Bbb C}}
\redefine\O{{\Cal O}}
\redefine\P{{\Bbb P}}
\redefine\phi{{\varphi}}
\redefine\epsilon{{\varepsilon}}
\define\Ker{\operatorname{Ker}}
\redefine\Im{\operatorname{Im}}
\define\ac{\acuteaccent}
\define\gr{\graveaccent}
\define\Res{\operatorname{Res}}
\define\set#1{\{#1\}}
\redefine\cdot{{\boldsymbol\cdot}}

\refstyle{B}
\NoRunningHeads
\TagsOnRight

\topmatter 
\title The third Cauchy-Fantappi\gr e formula of Leray \endtitle
\author Cezar Joi\c ta and Finnur L\ac arusson \endauthor
\affil University of Western Ontario \endaffil
\address Department of Mathematics, University of Western Ontario,
London, Ontario N6A~5B7, Canada \endaddress
\email cjoita\@uwo.ca, larusson\@uwo.ca \endemail
\date 3 January 2002 \enddate 

\thanks The second-named author was supported in part by the Natural
Sciences and Engineering Research Council of Canada.  \endthanks

\abstract We study the third Cauchy-Fantappi\gr e formula, an integral
representation formula for holomorphic functions on a domain in affine
space, presented by Jean Leray in the third paper of his famous series
{\it Probl\gr eme de Cauchy}, published in 1959.  We show by means of
examples that this formula does not hold without some additional
conditions, left unmentioned by Leray.  We give sufficient conditions
and a necessary condition for the formula to hold, and, in the case of a
contractible domain, characterize it cohomologically.  \endabstract


\endtopmatter

\document

\specialhead Introduction \endspecialhead

\noindent In the third paper \cite{L1} in his famous series {\it
Probl\gr eme de Cauchy}, Jean Leray founded the modern theory of
residues.  This paper is reprinted in Leray's {\it \OE uvres
scientifiques} \cite{L2} with an introduction by G\.  M\.  Henkin.  In
it, almost as an aside, Leray presents three representation formulas for
holomorphic functions on a domain in affine space, and calls them the
first, second, and third Cauchy-Fantappi\gr e formulas.  The first
formula, nowadays often called the Cauchy-Leray formula, is truly
fundamental: most integral representation formulas can be derived from
it in one way or another.  For an introduction to this area of complex
analysis, see \cite{B} and \cite{K}.  The second and third formulas,
obtained from the first one using residue theory, have received little
attention in the literature.  The proof of the second formula is
straightforward, but the third one turns out to be quite subtle.  We
will show by means of examples that it does not hold without some
additional conditions, left unmentioned by Leray.  We present and study
sufficient conditions and a necessary condition for the third formula to
hold, and, in the case of a contractible domain, characterize it
cohomologically. 

We assume that the reader is familiar with approximately the first half
of Leray's paper \cite{L1}, including the coboundary map and the
associated long exact homology sequence, absolute and relative residues,
the residue formula, and the interaction of the residue map with several
natural cohomology maps.  To establish a context and notation, we begin
by reviewing Leray's derivation of the three formulas.

\specialhead 1. The Cauchy-Fantappi\gr e formulas \endspecialhead

\noindent We start with a number of definitions, following Leray.  Let
$X$ be a domain in $\C^n$, $n\geq 1$, and let $Y=\P^n\times X$.  Leray
assumes that $X$ is convex, but we do not.  We define
$$Q=\set{(\xi,x)\in Y:\xi\cdot x = \xi_0+\xi_1 x_1+\dots+\xi_n x_n =
0}$$
and
$$P_z=\set{(\xi,x)\in Y:\xi\cdot z=0} = \text{hyperplane}\times X$$
for $z\in X$, which will be thought of as a fixed base point throughout. 
These are smooth hypersurfaces in $Y$ intersecting
transversely.  Here, $\xi_0,\dots,\xi_n$ are homogeneous coordinates
for a point $\xi\in\P^n$.  Define also
$$\omega(x)=dx_1\wedge\dots\wedge dx_n,$$
$$\omega'(\xi)=\sum_{k=1}^n (-1)^{k-1}\xi_k d\xi_1\wedge \dots \wedge
\widehat{d\xi_k}\wedge \dots \wedge d\xi_n,$$
and
$$\omega^*(\xi)=\sum_{k=0}^n (-1)^k\xi_k d\xi_0 \wedge \dots \wedge
\widehat{d\xi_k}\wedge \dots \wedge d\xi_n.$$
Then $(\xi\cdot z)^{-n}\omega'(\xi)$ is a well defined holomorphic
$(n-1)$-form and $(\xi\cdot z)^{-n-1}\omega^*(\xi)$ is a well defined
holomorphic $n$-form on $\P^n\setminus\set{\xi:\xi\cdot z=0}$.  Let
$$\Phi_z=\frac {\omega'(\xi)\wedge\omega(x)}{(\xi\cdot z)^n} \qquad
\text{and} \qquad \Psi_z=\frac {\omega^*(\xi)\wedge\omega(x)}
{(\xi\cdot z)^{n+1}}.$$
Then $\Phi_z$ is a holomorphic $(2n-1)$-form and $\Psi_z$ is a holomorphic
$2n$-form on $Y\setminus P_z$, with poles of order $n$ and $n+1$ along
$P_z$ respectively, and 
$$d\Phi_z=n\Psi_z.$$
We usually omit the subscript $z$.

We will work with homology and cohomology groups throughout the paper. 
Leray uses integral homology with the torsion part removed; we might as
well use homology with complex coefficients.  Our cohomology groups will
be complex de Rham groups, both absolute and relative. 

The projection $Q\setminus P\to X\setminus\set{z}$ is a fibration with
contractible fibre $\C^{n-1}$, so it is a homotopy equivalence and
induces isomorphisms on all homology and cohomology groups.  The
boundary map $H_{2n}(X,X\setminus\set{z}) \to
H_{2n-1}(X\setminus\set{z})$ is injective with a one-dimensional source. 
Its image, i.e., the kernel of the map $H_{2n-1}(X\setminus\set{z}) \to
H_{2n-1}(X)$, is generated by the image of a unique integral class
$\alpha_z$ in $H_{2n-1}(Q\setminus P)$ on which $i^{-n}\Phi>0$.  The
image of $\alpha$ in $H_{2n-1}(X\setminus\set{z})$ is represented by the
boundary of any smoothly bounded relatively compact subdomain of $X$
containing $z$, with the appropriate orientation.  It may be checked
that the orientation is given by the outward normal when $n\equiv 1,2
\pmod{4}$ and by the inward normal when $n\equiv 0,3 \pmod{4}$. 

With these definitions we can now state the first formula.  For the
proof, see \cite{L1, \S 56}.

\proclaim{The first Cauchy-Fantappi\gr e formula}  If $f$ is a
holomorphic function on a domain $X$ in $\C^n$ and $z\in X$, then
$$f(z) = \frac {(n-1)!}{(2\pi i)^n} \int_{\alpha_z} f\Phi_z.$$
\endproclaim

Here, by $f\Phi_z$, we mean $\Phi$ multiplied by $f$ precomposed with
the projection $Y\to X$.  Note that $f\Phi$ is closed on $Q\setminus P$,
being holomorphic of degree equal to the dimension of $Q\setminus P$, so
the integral makes sense. 

For $n=1$, this is the Cauchy formula.  In most applications, $\alpha$
is represented by a cycle which is a graph over the boundary of a
subdomain containing $z$.  This case is well described in Berndtsson's
survey \cite{B}, as are variants with weight factors and related
solution formulas for the $\overline\partial$ operator.  For an example
where $\alpha$ is represented by a cycle which is not a graph, see
\cite{K, 4.5}. 

Now we move on to the second formula.  We wish to express $\alpha$ as
the image of a class in $H_{2n-2}(Q\cap P)$ by Leray's coboundary
map $\delta: H_{2n-2}(Q\cap P)\to H_{2n-1}(Q\setminus P)$.  Let us
represent $\alpha$ by the $(2n-1)$-real-dimensional smooth submanifold
$$M=\set{([x\cdot\bar z-|x|^2,\bar x_1-\bar z_1,\dots, \bar x_n-\bar
z_n],x) : x\in\partial B}$$ 
of $Y$ with the orientation specified above, where $B\Subset X$ is the
open ball of radius $\epsilon>0$ centred at $z$.  Note that $\xi\cdot
x=0$ and $\xi\cdot z = -\epsilon^2$ for $(\xi,x)\in M$, so $M\subset
Q\setminus P$.  Using the same formula on $B\setminus\set{z}$ and taking
closure in $Q$, we get a $2n$-real-dimensional smooth submanifold of
$Q$, isomorphic to $\overline B$ with $z$ blown up, with boundary $M$. 
It intersects $P$ transversely in the common fibre $\P^{n-1}$ of $P$ and
$Q$ over $z$.  This fibre represents a class $\beta_z\in H_{2n-2}(Q\cap
P)$, oriented so that $\delta\beta = -\alpha$.  Applying the residue
formula (which says that the residue map is the dual of $\delta/2\pi i$)
to the first formula, we obtain immediately the second formula. 

\proclaim{The second Cauchy-Fantappi\gr e formula} If $f$ is a
holomorphic function on a domain $X$ in $\C^n$ and $z\in X$, then
$$f(z) = -\frac {(n-1)!}{(2\pi i)^{n-1}} \int_{\beta_z} \Res_{Q\cap P}
(f\Phi_z).$$
\endproclaim

By Leray's long exact sequence, the kernel of $\delta$ is the image of
the map $H_{2n}(Q)\to H_{2n-2}(Q\cap P)$ induced by intersecting with
$P$, which is injective since $H_{2n}(Q\setminus P)=0$.  Hence, $\beta$
is the unique class with $\delta\beta=-\alpha$ if and only if
$H_{2n}(Q)=0$, which is the case if $X$ is contractible. 

The class $\beta$ is never a boundary in $P$.  To obtain the third
formula we must express $\beta$ as a relative boundary with respect to a
new smooth hypersurface $S$ in $Y$, such that $P$, $Q$, and $S$ are in
general position.  Consider the following diagram.  
$$ \CD @.  @. 
H_{2n-2}(P\cap S) \\ @.  @.  @VVV \\ @.  \beta\in H_{2n-2}(Q\cap P) @>>>
H_{2n-2}(P) \\ @.  @VpVV @VVV \\ H_{2n-1}(P, Q\cup S) @>\partial>>
H_{2n-2}(Q\cap P, S) @>>> H_{2n-2}(P,S) \endCD $$ 
Assume that $p\beta=\partial\gamma$ for some $\gamma\in H_{2n-1}(P,Q\cup
S)$, i.e., that $\beta$ is homologous in $P$ to a cycle in $P\cap S$. 
We will refer to this as the {\it topological condition} on $S$.  It
holds if the map $H_{2n-2}(P\cap S)\to H_{2n-2}(P)$ is surjective.  It
also holds if, but not only if, $S$ contains a fibre of the projection 
$P\to X$: in Example C of Section 3, the inclusion $P\cap S
\hookrightarrow P$ is a homotopy equivalence, but $S$ does not contain 
a $P$-fibre.  Note that $\gamma$ is uniquely determined modulo the image
of $H_{2n-1}(P,S)$. 

Leray shows that the residue map $\Res$ commutes with $p^*$ and
anti-commutes with the dual $\partial^*$ of the boundary map.  Note that
$f\Phi$ vanishes on $(Q\setminus P)\cap S$ for dimensional reasons, so
it represents a class in $H^{2n-1}(Q\setminus P,S)$.  Now by the second
formula, 
$$\align
-\frac {(2\pi i)^{n-1}}{(n-1)!} f(z) &= \int_\beta\Res_{Q\cap P}f\Phi =
\int_\beta \Res_{Q\cap P}p^*(f\Phi) = \int_\beta p^*\Res_{(Q\cap P,S)} f\Phi
\\ 
&= \int_{p\beta} \Res_{(Q\cap P,S)} f\Phi = \int_{\partial\gamma} 
\Res_{(Q\cap P,S)}f\Phi 
= \int_\gamma \partial^*\Res_{(Q\cap P,S)}f\Phi \\
&= -\int_\gamma \Res_{(P,Q\cup S)}\partial^*(f\Phi).
\endalign $$
The map $\partial^*$ is the dual of the boundary map $\partial$ in {\it
relative} homology.  It is not simply the exterior derivative.  To
calculate $\partial^*(f\Phi)\in H^{2n}(Y\setminus P,Q\cup S)$, we extend
$f\Phi$ to a smooth form on $Y\setminus P$, vanishing on $S\setminus P$,
and differentiate.  If $\Phi|S\setminus P=0$, then the extension can be
chosen to be $f\Phi$ itself, and $\partial^*(f\Phi)$ is represented by
$d(f\Phi)=fd\Phi=nf\Psi$.  On page 155 of \cite{L1}, Leray seems to
assume that this is so, $\partial^*$ becomes $d$, and his third formula
reads as follows.

\proclaim{The third Cauchy-Fantappi\gr e formula} Let $f$ be a
holomorphic function on a domain $X$ in $\C^n$ and $z\in X$.  Let $S$ be
a smooth hypersurface in $Y$ such that $P_z$, $Q$, and $S$ are in
general position.  Then
$$f(z)=\frac{n!}{(2\pi i)^{n-1}}\int_\gamma\Res_{(P,Q\cup S)} f\Psi_z$$
for every $\gamma\in\partial^{-1}(p\beta_z)\subset H_{2n-1}(P,Q\cup S)$.
\endproclaim

Here, we view $f\Psi$ as a class in $H^{2n}(Y\setminus P,Q\cup S)$,
which the residue map takes to $H^{2n-1}(P,Q\cup S)$.  At this point,
the reader might find it helpful to consult Example A in Section 3 and
see what the formula looks like in a simple case.

Without some additional conditions on $S$, {\bf the third formula is
false}, as shown by Examples A, D, and E in Section 3.  The third
formula fails in two different ways in these examples.  In Examples D
and E, $n=2$ and the integral on the right hand side is not independent
of the choice of $\gamma$ in $\partial^{-1}(p\beta)$ for $f=1$.  In
Example A, $X=\C$ and $H_1(P,Q\cup S)$ is one-dimensional, so there is
only one choice of $\gamma$, but the formula fails for $f\in\O(\C)$ with
$f(1)\neq 0$.  We do not have an example where the formula holds for
$f=1$ but fails for other $f$.  In fact, we know nothing about the set
of functions for which the formula holds beyond the obvious fact that it
is a vector subspace of $\O(X)$. 

An important special case in which the third formula does hold is when
$S$ is the preimage of a smooth hypersurface $S_0$ in $X$.  Then
$\Phi|S\setminus P$ clearly vanishes, the topological condition is
satisfied, and it is easy to see that $P_z$, $Q$, and $S$ are in general
position if and only if $z\notin S_0$.  Examples B and C in Section 3
show that the third formula holds not only when $S$ is a preimage. 

Before presenting our examples, we give various answers to the question
of when the third formula holds, starting with a sufficient condition
which is also necessary in the contractible case.

\specialhead 2. Necessary and sufficient conditions
\endspecialhead

\noindent
Throughout this section, $X$ is a domain in $\C^n$, $z\in X$,
$f\in\O(X)$, $S$ is a smooth hypersurface in $Y$ such that $P_z$, $Q$,
and $S$ are in general position, and $S$ satisfies the topological
condition $\partial^{-1}(p\beta_z)\neq\varnothing$.

From our discussion in the previous section, it is clear that the third
formula holds if and only if $$\Res_{(P,Q\cup S)}\partial_Q^*(f\Phi) =
\Res_{(P,Q\cup S)} d(f\Phi) \qquad\text{on }
\partial^{-1}(p\beta)\subset H_{2n-1}(P,Q\cup S),$$ where we have
written $\partial_Q$ for the boundary map $H_{2n}(Y\setminus P,Q\cup
S)\to H_{2n-1}(Q\setminus P,S)$ previously simply denoted $\partial$,
because now we also want to consider the boundary map
$\partial_S:H_{2n}(Y\setminus P,Q\cup S)\to H_{2n-1}(S \setminus P,Q)$. 
For convenience, let us write $\phi$ for $f\Phi$.  As described above,
$\partial_Q^*\phi$ is represented by $d\sigma$, where $\sigma|Q\setminus
P=\phi$ and $\sigma|S\setminus P=0$.  Now $\phi$ also represents a class
in $H^{2n-1}(S\setminus P,Q)$, and $\partial_S^*\phi$ is represented by
$d\tau$, where $\tau|Q\setminus P=0$ and $\tau|S\setminus P=\phi$. 
Hence, $d\phi-\partial_Q^*\phi=\partial_S^*\phi$. 

Consider the following diagram.
$$\CD
@. H^{2n-2}(P,Q) @>{p^*}>> H^{2n-2}(P) \\
@. @V{i^*}VV @V{i^*}VV \\
H^{2n-3}(P\cap Q\cap S) @>{\partial^*}>> H^{2n-2}(P\cap S,Q) @>{p^*}>>
H^{2n-2}(P\cap S) \\
@V{\partial^*}VV  @V{\partial^*}VV  @. \\
H^{2n-2}(P\cap Q,S) @>{\partial^*}>> H^{2n-1}(P,Q\cup S) @. \\
@V{p^*}VV  @V{p^*}VV  @. \\
H^{2n-2}(P\cap Q) @>{\partial^*}>> H^{2n-1}(P,Q)  @.
\endCD $$
The rows and columns are parts of long exact sequences.  The top and
bottom squares commute.  Let us verify that the middle square
anti-commutes, pretending for a moment, to simplify the notation, that
$Q$ and $S$ are hypersurfaces in $P$.  Take a closed form $\omega_0$ on
$Q\cap S$ and extend it to a form $\omega$ on $P$.  The upper part of
the square works like this: first, $\partial^*\omega_0=d\omega|S$; then
we extend $d\omega|S$ to a form $\omega_1$ on $P$ such that
$\omega_1|Q=0$, and $\partial^*(d\omega|S)=d\omega_1$.  The lower part
of the square works like this: first, $\partial^*\omega_0=d\omega|Q$;
then we extend $d\omega|Q$ to a form $\omega_2$ on $P$ such that
$\omega_2|S=0$, and $\partial^*(d\omega|Q)=d\omega_2$.  We must show
that $d\omega_1+d\omega_2=0$ in $H^{2n-1}(P,Q\cup S)$.  This holds since
$d\omega_1+d\omega_2 = d(\omega_1+\omega_2-d\omega)$, and
$\omega_1+\omega_2-d\omega=0$ on $Q\cup S$. 

Write $\rho=\Res_{(P,Q\cup S)}\partial_S^*\phi$.  We have seen that the
third formula holds if and only if $\rho=0$ on $\partial^{-1}(p\beta)$,
i.e., $\rho=0$ on $\gamma$ and on $\Ker\partial$, where $\gamma$ is some
element of $\partial^{-1}(p\beta)$.  This is equivalent to
$\rho\in\Im\partial^*$, say $\rho=\partial^*\sigma$ with $\sigma\in
H^{2n-2}(P\cap Q,S)$, and $p^*\sigma(\beta)=\rho(\gamma) =0$. 

Let us now assume that $X$ is contractible.  Then $H_{2n-2}(P\cap Q)$ is
generated by $\beta$, so $p^*\sigma(\beta)=0$ implies $p^*\sigma=0$. 
Hence, the third formula is equivalent to $\rho\in\partial^*(\Ker p^*)$,
i.e., to $\rho$ coming from $H^{2n-3}(P\cap Q\cap S)$ (which is of
course trivial when $n=1$).  

The anti-commuting square
$$\CD
\phi\in H^{2n-1}(S\setminus P,Q) @>{\partial_S^*}>> H^{2n}(Y\setminus
P, Q\cup S) \\
@V{\Res}VV  @V{\Res}VV \\
H^{2n-2}(P\cap S, Q) @>{\partial^*}>> H^{2n-1}(P,Q\cup S)\owns\rho
\endCD $$
shows that $\rho$ comes from $H^{2n-2}(P\cap S, Q)$ with preimage
$-\Res_{(P\cap S, Q)}\phi$ there.  By chasing the big diagram, we see
that $\rho$ comes from $H^{2n-3}(P\cap Q\cap S)$ if and only if
$p^*\Res_{(P\cap S,Q)}\phi = \Res_{P\cap S}\phi$ comes from
$H^{2n-2}(P,Q)$.  Now the map $p^*:H^{2n-2}(P,Q)\to H^{2n-2}(P)$ is
zero, because the map $i^*:H^{2n-2}(P)\to H^{2n-2}(P\cap Q)$, which
follows it in the long exact cohomology sequence, is an isomorphism. 
Hence, the third formula holds if and only if $\Res_{P\cap S}\phi$ is
zero in $H^{2n-2}(P\cap S)$ (here, we view $\phi$ as representing a
class in $H^{2n-1}(S\setminus P)$).  This is equivalent to the existence
of a smooth $(2n-2)$-form $\sigma$ on $S\setminus P$ such that
$\phi+d\sigma$ extends smoothly to $S$.  In the general case, when $X$
is not necessarily contractible, this is a sufficient condition for the
formula to hold. 

Let us now derive a necessary condition for the third formula to hold. 
This condition fails in Examples D and E in Section 3.  The third
formula implies that the integral on its right hand side is independent
of the choice of $\gamma$ in $\partial^{-1}(p\beta)$, i.e., that
$$\int_{p\theta} \Res_{(P,Q\cup S)}d\phi=0 \qquad \text{for every }
\theta\in H_{2n-1}(P,S).$$
Viewing $d\phi$ as also representing a class in $H^{2n}(Y\setminus P,S)$
that is the image by $p^*$ of the class represented by $d\phi$ in
$H^{2n}(Y\setminus P, Q\cup S)$, we see that this condition is
equivalent to $\Res_{(P,S)}d\phi$ vanishing in $H^{2n-1}(P,S)$.  This is
equivalent to $d\phi$ coming from $H^{2n}(Y,S)$, which means that there
is a smooth $(2n-1)$-form $\sigma$ on $Y\setminus P$ with
$\sigma|S\setminus P=0$ such that $d(\phi-\sigma)$ extends smoothly to
$Y$. 

We have proved the following result.

\proclaim{Theorem}  Let $X$ be a domain in $\C^n$, $z\in X$,
$f\in\O(X)$, and $S$ be a smooth hypersurface in $Y$ such that $P_z$, $Q$,
and $S$ are in general position and $S$ satisfies the topological
condition.  For Leray's third formula to hold for this data, it is
necessary that 
$$\Res_{(P_z,S)}f\Psi_z=0 \qquad \text{in } H^{2n-1}(P_z,S),$$
and sufficient that
$$\Res_{P_z\cap S}f\Phi_z=0 \qquad \text{in } H^{2n-2}(P_z\cap S).$$
When $X$ is contractible, this sufficient condition is also necessary.
\endproclaim

In order to get a characterization of the third formula, we do not
actually need $X$ to be contractible.  The argument above shows that it
suffices to have $H_{2n-2}(P\cap Q)$ and
$H_{2n-2}(P)=H_{2n-2}(\P^{n-1}\times X)$ one-dimensional (each generated
by $\beta$). 

There is no shortage of other sufficient conditions for Leray's third
formula to hold.  One is that $\phi|S\setminus P$ extend holomorphically
to $S$: then, surely, its residue along $P\cap S$ vanishes.  Examples B
and C in Section 3 illustrate this.  Another is for $H_{2n-1}(S\setminus
P)\to H_{2n-1}(S)$ to be injective: then every class in
$H^{2n-1}(S\setminus P)$ lies in the image of $H^{2n-1}(S)$.  A third
sufficient condition is that $\partial_Q^*\phi$ and $d\phi$ have the
same residue along $(P,Q\cup S)$, as is evident from our discussion
preceding the statement of the formula.  This is equivalent to the
existence of a smooth $(2n-1)$-form $\sigma$ on $Y\setminus P$ with
$\sigma|S\setminus P=0$ and $\sigma=\phi$ on $Q\setminus P$ such that
$d(\phi-\sigma)$ extends smoothly to $Y$. 

We have pinched Leray's third formula between this sufficient condition
and the necessary condition that such a form $\sigma$ exist without
requiring that $\sigma=\phi$ on $Q\setminus P$.  We have also pinched
Leray's third formula between the sufficient condition that $\Res_{P\cap
S}\phi=0$ and the (same) necessary condition that $\partial^*\Res_{P\cap
S}\phi = -\Res_{(P,S)}d\phi=0$. 

We conclude this section by pointing out that if $S$ contains a fibre of
the projection $P\to X$ (which by Example C of Section 3 is stronger
than the topological condition), then the sufficient condition
$\Phi|S\setminus P=0$ actually implies our strongest and most obvious
sufficient condition, that $S$ be the preimage of a hypersurface in $X$. 
We may assume that $z=0$.  First, in the affine coordinates
$y_j=\xi_j/\xi_0$, we have
$$\Phi=y_1^n d(\frac{y_2}{y_1})\wedge\dots\wedge d(\frac{y_n}{y_1})\wedge
\omega(x),$$
so, changing to the coordinates $w_1=y_1$ and $w_j=y_j/y_1$ for
$j=2,\dots,n$, we have $\Phi=0$ on the subset of $S\setminus P$ where
$\xi_1\neq 0$ if and only if $dw_2\wedge\dots\wedge dw_n\wedge
dx_1\wedge\dots\wedge dx_n=0$ there.  This is equivalent to $\partial
s/\partial w_1=0$ for any local defining function $s$ for $S$, i.e., to
$S$ consisting of curves on which $w_2,\dots,w_n,x_1,\dots,x_n$ are
constant.  Hence, for every $x\in X$, $S\cap(\C^n\times\set{x})$ is a
union of lines through the origin, or the origin itself, or empty. 
Adding the point at infinity to each of these affine lines yields the
closure of $S\setminus P$ in $Y$, which equals $S$ since $S\setminus P$
is dense in $S$ (this clearly follows from $P$ and $S$ being in general
position, but $P\not\subset S$ is enough).  Hence, the assumption that
$S$ contains a $P$-fibre implies that $S$ contains the whole fibre
$\P^n\times\set{x}$ for some $x\in X$.  Being a hypersurface, $S$ cannot
also contain $\set{0}\times X$, so its image by the projection $Y\to X$
(which is proper and hence takes a subvariety to a subvariety) must be a
hypersurface in $X$ with $S$ as its preimage in $Y$.

\specialhead 3.  Examples \endspecialhead

\noindent
In this section, we present the following five examples of smooth
hypersurfaces $S$ in $\P^n\times X$ such that $P_z$, $Q$, and $S$ are in
general position for $z=0$, and $S$ contains a fibre of the projection
$P\to X$, so $S$ satisfies the topological condition, i.e., the basic
assumptions in the third Cauchy-Fantappi\gr e formula are in place. 
\newline 
{\bf A.}  Examples with $X=\C$ where the third formula holds (to show what
it looks like in the simplest case) and where it fails. \newline
{\bf B.}  A more involved example with $X=\C$ where $\Phi|S\setminus P$ 
extends holomorphically to $S$, so the third formula holds, although $S$
is not the preimage of a subset of $X$. \newline
{\bf C.}  Same as Example B, but with $X=\C^2$.  Also, a modification of
this example that satisfies the topological condition even though it
does not contain a $P$-fibre.  \newline
{\bf D.}  An example with $n=2$ where the necessary condition
$\Res_{(P,S)}f\Psi=0$, and hence the third formula, fails for $f=1$; 
here, $X$ is contractible but not convex. \newline
{\bf E.} An example with $X=\C^2$ where the necessary condition fails
for $f=1$; here, the computations are considerably more complicated than
in Example D.

\smallskip\noindent {\bf Example A.} We first consider the third
Cauchy-Fantappi\gr e formula in the simplest case, with $X=\C$ and
$z=0$.  Letting $\eta=\xi_0/\xi_1$, we have $P=\set{\eta=0}\cong\C$,
$Q=\set{\eta+x=0}$, and $\Psi=\eta^{-2} dx\wedge d\eta$.  Take
$S=\set{a\eta+x=1}$ with $a\in\C$.  Then $P$, $Q$, and $S$ are in
general position and $S$ satisfies the topological condition.  In fact,
$S$ contains the $P$-fibre above $x=1$, and $\gamma\in H_1(P,\set{0,1})$
is represented by any path $\langle 1,0\rangle$ in $\C$ from $1$ to 
$0$.  Let
$$\align
\sigma&=\frac{a\eta+x-1}\eta(d\eta+dx)-\frac{\eta+x}\eta(ad\eta+dx) \\  
&=((1-a)\frac x \eta - \frac 1 \eta)d\eta + (a-1-\frac 1 \eta)dx.
\endalign$$
The first expression for $\sigma$ shows that it vanishes on $Q\cup S$,
and the second one shows that 
$$d\sigma=(a-1)\frac{d\eta}\eta\wedge dx -\Psi,$$
so $\Res_{(P,Q\cup S)}\Psi$ is represented by the form $(a-1)dx$.
If $f\in\O(\C)$, then similarly $f\sigma$ vanishes on $Q\cup S$ and
$$f\Psi + d(f\sigma)=\frac{d\eta}\eta\wedge d(f(x)((a-1)x+1)),$$
so $\Res_{(P,Q\cup S)}f\Psi$ is represented by the form
$d(f(x)((a-1)x+1))$, and the third formula looks as follows:
$$f(0)=\int_{\langle 1,0\rangle} d(f(x)((a-1)x+1))
= f(0) - af(1).$$
Clearly, the formula holds for all $f\in\O(\C)$ only when $a=0$, but
fails for $f$ with $f(1)\neq 0$ when $a\neq 0$.

\smallskip\noindent{\bf Example B.}  Again, we take $X=\C$ and $z=0$,
but now $S$ is defined by the equation
$$s(\xi_0,\xi_1,x)=\xi_0^2+\xi_1(\xi_0+\xi_1)(x-1)=0.$$
Clearly, $S$ is not the preimage of a subset of $\C$, and $P\cap S$ is 
the point above $1$ in $\C$.  It is easily verified that $S$ is smooth
and $P$, $Q$, and $S$ are in general position.  

We claim that $\Phi|S\setminus P$, and hence $f\Phi|S\setminus P$ for
every $f\in\O(X)$, extends holomorphically to $S$.  Working near $P$,
where $\xi_1\neq 0$, and using the affine coordinate $\eta=\xi_0/\xi_1$,
we have $\Phi=dx/\eta$.  On $S$ near $P$, $x=1-\eta^2/(\eta+1)$, so
$dx=-\eta(\eta+2)/(\eta+1)^2 d\eta$, and $\Phi=-(\eta+2)/(\eta+1)^2 d\eta$,
which extends holomorphically across $\eta=0$.  Hence, the third formula
holds.  Let us check this by calculating the integrand $\Res_{(P,Q\cup
S)}f\Psi$.  We work on the neighbourhood of $P$ where $\xi_1\neq 0$. 
Let $q$ be the defining function $\eta + x$ for $Q$.  Let
$$\sigma = \frac f \eta(s dq - q ds) = -\frac{f(x)(x-1)^2}\eta
d\eta - \frac f \eta dx + \theta,$$
where $\theta$ is smooth near $P$.  Then $\sigma$ vanishes on $Q\cup S$ and
$$f\Psi+d\sigma = \frac{d\eta}\eta\wedge d(f(x)(x-1)^2)+d\theta,$$
so the residue $\Res_{(P,Q\cup S)}f\Psi$ is represented by the form
$d(f(x)(x-1)^2)$.  Integrating this form along a path from $1$ to $0$ gives
$f(0)$. 

\smallskip\noindent{\bf Example C.}  Here, $X=\C^2$, $z=0$, and $S$ is 
defined by the equation
$$s(\xi_0,\xi_1,\xi_2,x_1,x_2)=\xi_0^3+\xi_1^3(x_1-1)+\xi_2^3(x_2-2)=0.$$
Note that $S$ is not the preimage of a curve in $X$, and $S$ contains the
$P$-fibre above $(1,2)$.  Lengthy but routine computations show that $S$
is smooth and $P$, $Q$, and $S$ are in general position.  We will show
that $\Phi|S\setminus P$ extends holomorphically to $S$, so the third 
formula holds.  Let $U_k$ be the subset of $S$ where
$\xi_k\neq 0$, $k=1,2$.  We shall verify that $\Phi|S\setminus P$
extends across $P$ on $U_2$; the case of $U_1$ is analogous.  In the
affine coordinates $y_0=\xi_0/\xi_2$, $y_1=\xi_1/\xi_2$, we have
$\Phi=-y_0^{-2} dy_1\wedge dx_1\wedge dx_2$.  On $S\cap U_2$,
$x_2=2-y_0^3-y_1^3(x_1-1)$, so there, $\Phi=3dy_0\wedge dy_1\wedge
dx_1$, which clearly extends across $P$.

Now take
$$s(\xi_0,\xi_1,\xi_2,x_1,x_2)=\xi_0^3+\xi_1^3(x_1-1)+\xi_2^3(x_2-2) 
+2\xi_1^2\xi_2.$$
Just as before, $S$ is smooth, $P$, $Q$, and $S$ are in general
position, and $\Phi|S\setminus P$ extends holomorphically to $S$.  Also,
$S$ is not the preimage of a curve in $X$: in fact, it is easy to see
that $S$ does not even contain a $P$-fibre.  All the same, $S$ satisfies
the topological condition, as we will now verify.  Consider the
projection $\pi:P\cap S \hookrightarrow P = \P^1\times X \to \P^1$. 
Since $P\to \P^1$ is a homotopy equivalence, if we can show that $\pi$ 
is a homotopy equivalence, then so is the inclusion $P\cap S
\hookrightarrow P$, and the topological condition follows.

We will show that $\pi$ is a fibration with contractible fibres.  Now
$P\cap S$ is given by the equations $\xi_0=0$ and $\xi_1^3(x_1-1)
+\xi_2^3(x_2-2)+2\xi_1^2\xi_2=0$.  The fibres of $\pi$ are lines in
$\C^2$ since the coefficients of $x_1$ and $x_2$ do not vanish
simultaneously.  We claim that $\pi$ is surjective.  Namely, if
$\xi=[\xi_1,\xi_2]\in\P^1$ and $\xi_1\neq 0$, we set $x_2=0$ and
$x_1=\xi_1^{-3}(\xi_1^3+2\xi_2^3-2\xi_1^2\xi_2)$.  Then $(\xi,x)\in
P\cap S$ and $\pi(\xi,x)=\xi$.  If $\xi_2\neq 0$, we take $x_1=0$ and
$x_2=\xi_2^{-3}(\xi_1^3+2\xi_2^3-2\xi_1^2\xi_2)$.  Finally, $\pi$ is
locally trivial: over $U=\{\xi\in\P^1:\xi_1\neq 0\}$, the trivialization
map is $\pi^{-1}(U) \to U\times\C$, $(\xi,x_1,x_2)\mapsto (\xi,x_2)$,
and the inverse map is given by $x_1=1-
\xi_1^{-3}(\xi_2^3(x_2-2)+2\xi_1^2\xi_2)$.  Over $\{\xi_2\neq 0\}$, the
trivialization is given by $(\xi,x_1,x_2)\mapsto (\xi,x_1)$ and the
inverse map is obtained in the same way. 

\smallskip\noindent {\bf Example D.}  Here, $n=2$, $z=0$, and $S$ is 
defined by the equation
$$s(\xi_0,\xi_1,\xi_2,x_1,x_2)=\xi_0^2+\xi_1(\xi_1+\xi_2)(x_1-1)x_2
+\xi_2^2(x_2^2+1) = 0$$
over a contractible domain $X$ in $\C^2$ which will be specified later. 
Note that $S$ contains the $P$-fibres over $(1,\pm i)$.  Computations,
whose details will not be reproduced here, show that $S$ is smooth
outside $P$; at each point of $P\cap S$, $S$ is smooth and in general
position with respect to $P$ except over the point $(1,0)$; at each
point of $(Q\cap S)\setminus P$, $Q$ and $S$ are in general position
except over a finite set in $\C^2$ which does not contain $(1,\pm i)$ or
the base point $(0,0)$; and at each point of $P\cap Q\cap S$, the three
surfaces are in general position. 

Let us show that the necessary condition fails, i.e.,
$\Res_{(P,S)}\Psi\neq 0$.  Let $U=\set{\xi_2\neq 0}$.  It suffices to
prove that $\Res_{(P\cap U,S\cap U)}\Psi|U \neq 0$.  In the following,
we will work on $U$ but omit it from the notation.  In the affine
coordinates $y_0=\xi_0/\xi_2$, $y_1=\xi_1/\xi_2$,
$$\Psi=y_0^{-3} dy_0\wedge dy_1\wedge dx_1\wedge dx_2$$
and
$$s=y_0^2 + y_1(y_1+1)(x_1-1)x_2 + x_2^2 +1.$$
Let 
$$\tau=\frac s {y_0^2} dy_1\wedge dx_1\wedge dx_2 + \frac
{(x_1-1)dy_1\wedge dx_2 - x_2 dy_1\wedge dx_1}{2y_0^2}\wedge ds.$$
Then $\tau|S$ clearly vanishes and
$$d\tau = 2(\frac 1{y_0}-\frac 1{y_0^3})dy_0\wedge dy_1\wedge dx_1\wedge
dx_2,$$
so 
$$\Psi + \frac 1 2 d\tau = \frac {dy_0} {y_0}\wedge dy_1\wedge dx_1\wedge
dx_2, $$
and $\Res_{(P,S)}\Psi$ is represented by the form $dy_1\wedge dx_1\wedge
dx_2 = d((1-x_1)dy_1\wedge dx_2)$.  We need to show that this form is not
nullcohomologous on $P$ relative to $S$.  By Leray's long exact 
sequence and contractibility of $X$, it suffices to show that
$(1-x_1)dy_1\wedge dx_2$ is not exact on $S\cap P$.  Now on $S\cap P$,
$$1-x_1 = \frac {x_2^2+1}{x_2 y_1(y_1+1)}$$
(where the denominator does not vanish),
and integrating the form $(1-x_1)dy_1\wedge dx_2$ over the two-cycle
$$y_1=\epsilon e^{i\theta}, \quad x_2=\epsilon e^{i\eta}, \quad
x_1 = 1 - \frac {1+\epsilon^2 e^{2i\eta}}{\epsilon^2 e^{i(\theta+\eta)}
(1+\epsilon e^{i\theta})}, \qquad \theta,\eta\in[0,2\pi],$$
in $S\cap P$ with $0<\epsilon<1$, using the one-dimensional residue
formula, we do not get zero.

So $X$ must satisfy the following conditions: $X$ is contractible, $X$
contains $(0,0)$ and $(1,\pm i)$, $X$ avoids a certain finite set not
containing these points, and $X$ contains the set of points
$(x_1,x_2)$ as above with $\theta,\eta\in[0,2\pi]$ for some
$\epsilon\in(0,1)$.  This set is a two-dimensional real submanifold of
$\C^2$ which does not disconnect $\C^2$.  We can for instance take $X$ to
be the complement of a broken half-line joining the points we must avoid
and going out to infinity.

\smallskip\noindent {\bf Example E.}  Here, $X=\C^2$, $z=0$, and $S$ is 
defined by the equation
$$s(\xi_0,\xi_1,\xi_2,x_1,x_2)=\xi_0^2 + (\xi_1^2 + 3\xi_1\xi_2 x_2 +
2\xi_2^2 x_2^2)(x_1-1) +\xi_2^2(x_2^3+1)=0.$$
Clearly, $S$ contains the three $P$-fibres above $(1,\root 3\of{-1})$. 
Calculations, which are considerably more complicated than the
corresponding ones in Example D and for which one may want to use a
computer algebra system, show that $S$ is smooth and $P$, $Q$, and $S$ 
are in general position.

Again, working on $U=\set{\xi_2\neq 0}$, which will be omitted from the
notation, with affine coordinates $y_0=\xi_0/\xi_2$, $y_1=\xi_1/\xi_2$,
we shall show that $\Res_{(P,S)}\Psi \neq 0$.  Let
$$\tau=\frac s {y_0^2} dy_1\wedge dx_1\wedge dx_2 - \frac {y_1
dx_1\wedge dx_2-(x_1-1)dy_1\wedge dx_2 + x_2 dy_1\wedge dx_1}{3y_0^2}
\wedge ds.$$ 
We see that $\tau|S=0$ and compute that
$$d\tau = 2(\frac 1 {y_0} - \frac 1{y_0^3}) dy_0\wedge dy_1\wedge
dx_1\wedge dx_2,$$
so as before, $\Res_{(P,S)}\Psi$ is represented by the form $dy_1\wedge
dx_1\wedge dx_2$, and we must show that $\theta=(1-x_1)dy_1\wedge dx_2$
is not exact on $S\cap P$.  

The projection $(y_1,x_1,x_2)\mapsto (y_1,x_2)$ restricts to an isomorphism 
of $S\cap P\setminus\set{(y_1+x_2)(y_1+2x_2)=0}$ 
onto $\C^2\setminus
\set{(y_1+x_2)(y_1+2x_2)=0}$ with inverse given by
$$1-x_1=\frac{x_2^3+1}{(y_1+x_2)(y_1+2x_2)}.$$  
Set $u=y_1+x_2$ and $v=y_1+2x_2$, so $x_2=v-u$, $y_1=2u-v$, and
$dy_1\wedge dx_2=du\wedge dv$.  The map $(y_1,x_1,x_2)\mapsto (u,v)$ is
an isomorphism of $S\cap P\setminus\set{(y_1+x_2)(y_1+2x_2)=0}$ onto
$\C^2\setminus\set{uv=0}$, and the push-forward
$\dfrac{(v-u)^3+1}{uv}du\wedge dv$ of $\theta$ is not exact on the
image: just integrate it over the product of two small circles centred
at the origin, using the one-dimensional residue formula.

\enddocument